\documentclass[12pt]{amsart}
\usepackage{amssymb}
\usepackage{amsfonts}
\usepackage[pdfpagemode=none,pdfstartview=FitH]{hyperref}

\newtheorem{theorem}{Theorem}[section]
\theoremstyle{plain}

\newtheorem{lemma}[theorem]{Lemma}

\newtheorem{proposition}[theorem]{Proposition}
\newtheorem{remark}[theorem]{Remark}

\numberwithin{equation}{section}
\input{tcilatex}

\begin{document}
\title[The Symmetric C*-algebra for two unitaries]{The C*-algebra of
symmetric words in two universal unitaries}
\author{Man-Duen Choi}
\address{Department of Mathematics\\
University of Toronto}
\email{choi@math.toronto.edu}
\author{Fr\'{e}d\'{e}ric Latr\'{e}moli\`{e}re}
\address{Department of Mathematics\\
University of Toronto}
\email{frederic@math.toronto.edu}
\urladdr{http://www.math.toronto.edu/frederic}
\date{October 2006}

\begin{abstract}
We compute the $K$-theory of the C*-algebra of symmetric words in two
universal unitaries. This algebra is the fixed point C*-algebra for the
order-two automorphism of the full C*-algebra of the free group on two
generators which switches the generators. Our calculations relate the $K$
-theory of this C*-algebra to the $K$-theory of the associated
C*-crossed-product by $\mathbb{Z}_{2}$.
\end{abstract}

\maketitle

\section{Introduction}

This paper investigates an example of a C*-algebra of symmetric words in
noncommutative variables. Our specific interest is in the C*-algebra of
symmetric words in the two universal unitaries generating the full
C*-algebra $C^{\ast }\left( \mathbb{F}_{2}\right) $ of the free group on two
generators. Our main result is the computation of the $K$-theory \cite%
{Blackadar98}\ of this algebra.

The two canonical unitary generators of $C^{\ast }\left( \mathbb{F}
_{2}\right) $ are denoted by $U$ and $V$. The C*-algebra of symmetric words
in two universal unitaries $U$,\ $V$ is precisely defined as the fixed point
C*-algebra $C^{\ast }\left( \mathbb{F}_{2}\right) _{1}$ for the order-2
automorphism $\sigma $ which maps $U$ to $V$ and $V$ to $U$. Our strategy to
compute the $K$-theory of $C^{\ast }\left( \mathbb{F}_{2}\right) _{1}$
relies upon the work of Rieffel \cite[Proposition 3.4]{Rieffel80} about
Morita equivalence between fixed point C*-algebras and C*-crossed-products.

The first part of this paper describes the two algebras of interest: the
fixed point C*-algebra $C^{\ast }\left( \mathbb{F}_{2}\right) _{1}$ and the
crossed-product $C^{\ast }\left( \mathbb{F}_{2}\right) \rtimes _{\sigma } 
\mathbb{Z}_{2}$ \cite{Zeller-Meier68}\cite{Davidson}. We also exhibit an
ideal $\mathcal{J}$ in $C^{\ast }\left( \mathbb{F}_{2}\right) \rtimes
_{\sigma }\mathbb{Z}_{2}$ which is strongly Morita equivalent to $C^{\ast
}\left( \mathbb{F}_{2}\right) _{1}$ and can be easily described as the
kernel of an very simple *-morphism. This allows us to reduce the problem to
the calculation of the $K$-theory of $C^{\ast }\left( \mathbb{F}_{2}\right)
\rtimes _{\sigma }\mathbb{Z}_{2}$.

The second part of this paper starts by this calculation. We use a standard
result from Cuntz \cite{Cuntz82} to calculate the $K$-theory of $C^{\ast
}\left( \mathbb{F}_{2}\right) \rtimes _{\sigma }\mathbb{Z}_{2}$. From this,
we compute the $K$-theory of $C^{\ast }\left( \mathbb{F}_{2}\right) _{1}$,
using the ideal $\mathcal{J}$, using simply a six-terms exact sequence in $K$
-theory.

We conclude the paper by looking a little closer to the obstruction to the
existence of unitaries of nontrivial $K$-theory in $C^{\ast }\left( \mathbb{%
F }_{2}\right) _{1}$. This amounts to comparing the structure of the ideal $%
\mathcal{J}$ and an ideal in $C^{\ast }\left( \mathbb{F}_{2}\right) _{1}$
related to the same representation as $\mathcal{J}$.

We also should mention that, in principle, using the results in \cite%
{Latremoliere06}, one could derive information on the representation theory
of $C^{\ast }\left( \mathbb{F}_{2}\right) \rtimes _{\sigma }\mathbb{Z}_{2}$
from the representations of $C^{\ast }\left( \mathbb{F}_{2}\right) $.

\section{The Fixed Point C*-algebra}

Let $C^{\ast }\left( \mathbb{F}_{2}\right) =C^{\ast }(U,V)$ be the universal
C*-algebra generated by two unitaries $U$ and $V$. We consider the order-2
automorphism $\sigma $ of $C^{\ast }\left( \mathbb{F}_{2}\right) $ uniquely
defined by $\sigma (U)=V$ and $\sigma (V)=U$. These relations indeed define
an automorphism by universality of $C^{\ast }\left( \mathbb{F}_{2}\right) $.
Our main object of interest is the fixed point C*-algebra $C^{\ast }\left( 
\mathbb{F}_{2}\right) _{1}=\left\{ a\in C^{\ast }\left( \mathbb{F}
_{2}\right) :\sigma (a)=a\right\} $ which can be seen as the C*-algebra of
symmetric words in two universal unitaries. The C*-algebra $C^{\ast }\left( 
\mathbb{F}_{2}\right) _{1}$ is related \cite{Rieffel80} to the
C*-crossed-product $C^{\ast }\left( \mathbb{F}_{2}\right) \rtimes _{\sigma } 
\mathbb{Z}_{2}$, which we will consider in our calculations. Our objective
is to gain some understanding of the structure of the unitaries and
projections in $C^{\ast }\left( \mathbb{F}_{2}\right) _{1}$.

The first step in our work is to describe concretely the two C*-algebras $%
C^{\ast }\left( \mathbb{F}_{2}\right) _{1}$ and $C^{\ast }\left( \mathbb{F}
_{2}\right) \rtimes _{\sigma }\mathbb{Z}_{2}$. Using \cite[Proposition 3.4]%
{Rieffel80}, we also exhibit an ideal in the crossed-product $C^{\ast
}\left( \mathbb{F}_{2}\right) \rtimes _{\sigma }\mathbb{Z}_{2}$ which is
Morita equivalent to $C^{\ast }\left( \mathbb{F}_{2}\right) _{1}$.

The following easy lemma will be useful in our work:

\begin{lemma}
\label{split}Let $A$ be a unital C*-algebra and $\sigma $ an order-two
automorphism of $A$. The fixed-point C*-algebra of $A$ for $\sigma $ is the
set $A_{1}=\left\{ a+\sigma (a):a\in A\right\} $. Set $A_{-1}=\left\{
a-\sigma (a):a\in A\right\} $. Then $A=A_{1}+A_{-1}$ and $A_{1}\cap
A_{-1}=\left\{ 0\right\} $.
\end{lemma}

\begin{proof}
Let $\omega \in A$. Since $\sigma ^{2}=1$ we have $\sigma (\omega +\sigma
(\omega ))=\sigma (\omega )+\omega \in A_{1}$. On the other hand, if $a\in A$
then $a=\frac{1}{2}\left( a+\sigma (a)\right) +\frac{1}{2}\left( a-\sigma
(a)\right) $. Hence if $a\in A$ then $a-\sigma (a)=0$ and $a\in \left\{
\omega +\sigma (\omega ):\omega \in A\right\} $. Moreover this proves that $%
A=A_{1}+A_{-1}$. Last, if $a\in A_{1}\cap A_{-1}$ then $a=\sigma (a)=-\sigma
(a)=0$.
\end{proof}

We can use Lemma (\ref{split}) to obtain a more concrete description of the
fixed-point C*-algebra $C^{\ast }\left( \mathbb{F}_{2}\right) _{1}$ of
symmetric words in two unitaries:

\begin{theorem}
Let $\sigma $ be the automorphism of $C^{\ast }\left( \mathbb{F}_{2}\right)
=C^{\ast }(U,V)$ defined by $\sigma (U)=V$ and $\sigma (V)=U$. Then the
fixed point C*-algebra of $\sigma $ is $C^{\ast }\left( \mathbb{F}
_{2}\right) _{1}=C^{\ast }\left( U^{n}+V^{n}:n\in \mathbb{N}\right) $.
\end{theorem}

\begin{proof}
Obviously $\sigma (U^{n}+V^{n})=U^{n}+V^{n}$ for all $n\in \mathbb{Z}$.
Hence 
\begin{equation*}
C^{\ast }\left( \left\{ U^{n}+V^{n}:n\in \mathbb{Z}\right\} \right)
\subseteq C^{\ast }\left( \mathbb{F}_{2}\right) _{1}\text{.}
\end{equation*}%
Conversely, $C^{\ast }\left( \mathbb{F}_{2}\right) _{1}=\left\{ \omega
+\sigma (\omega ):\omega \in C^{\ast }\left( \mathbb{F}_{2}\right) \right\} $
by Lemma (\ref{split}). So $C^{\ast }\left( \mathbb{F}_{2}\right) _{1}$ is
generated by elements of the form $\omega +\sigma (\omega )$ where $\omega $
is a word in $C^{\ast }\left( \mathbb{F}_{2}\right) $, since $C^{\ast
}\left( \mathbb{F}_{2}\right) $ is generated by words in $U$ and $V$, i.e.
by monomial of the form$U^{a_{0}}V^{a_{1}}\ldots U^{a_{n}}$ with $n\in 
\mathbb{N}$, $a_{0},\ldots ,a_{n}\in \mathbb{Z}$. It is thus enough to show
that for any word $\omega \in C^{\ast }\left( \mathbb{F}_{2}\right) $ we
have $\omega +\sigma (\omega )\in S$ where $S=C^{\ast }\left(
U^{n}+V^{n}:n\in \mathbb{Z}\right) $. Since, if $\omega $ starts with a
power of $V$ then $\sigma (\omega )$ starts with a power of $U$, we may as
well assume that $\omega $ always starts with a power of $U$ by symmetry.
Since the result is trivial for $\omega =1$ we assume that $\omega $ starts
with a nontrivial power of $U$. Such a word is of the form $\omega
=U^{a_{0}}V^{a_{1}}\ldots U^{a_{n-1}}V^{a_{n}}$ with $a_{0},\ldots
,a_{n-1}\in \mathbb{Z}\backslash \left\{ 0\right\} $ and $a_{n}\in \mathbb{Z}
$.

We define the order of such a word $\omega $ as the integer $o(\omega )=n$
if $a_{n}\not=0$ and $o(\omega )=n-1$ otherwise. In other words, $o(\omega )$
is the number of times we go from $U$ to $V$ or $V$ to $U$ in $\omega $. The
proof of our result follows from the following induction on $o(\omega )$.

By definition, if $\omega $ is a word such that $o(\omega )=0$ then $\omega
+\sigma (\omega )\in S$. Let us now assume that for some $m\geq 1$ we have
shown that for all words $\omega $ starting in $U$ such that $o(\omega )\leq
m-1$ we have $\omega +\sigma (\omega )\in S$. Let $\omega $ be a word of
order $m$ and let us write $\omega =U^{a_{0}}\omega _{1}$ with $\omega _{1}$
a word starting in a power of $V$. By construction, $\omega _{1}$ is of
order $m-1$. Let $\omega _{2}=V^{a_{0}}\omega _{1}$. By construction, $%
o(\omega _{2})=m-1$ or $m-2$. Either way, by our induction hypothesis, we
have $\omega _{1}+\sigma (\omega _{1})\in S$ and $\omega _{2}+\sigma (\omega
_{2})\in S$. Now:%
\begin{eqnarray*}
\omega +\sigma (\omega ) &=&U^{a_{0}}\omega _{1}+V^{a_{0}}\sigma (\omega
_{1}) \\
&=&\left( U^{a_{0}}+V^{a_{0}}\right) \left( \omega _{1}+\sigma (\omega
_{1})\right) -\left( \omega _{2}+\sigma (\omega _{2})\right)
\end{eqnarray*}%
hence $\omega +\sigma (\omega )\in S$ and our induction is complete. Hence $%
C^{\ast }\left( \mathbb{F}_{2}\right) _{1}=S$ as desired.
\end{proof}

We wish to understand more of the structure of the fixed point C*-algebra $%
C^{\ast }\left( \mathbb{F}_{2}\right) _{1}$. Using Morita equivalence, we
can derive its $K$-theory. According to \cite[Proposition 3.4]{Rieffel80}, $%
C^{\ast }\left( \mathbb{F}_{2}\right) _{1}$ is strongly Morita equivalent to
the ideal $\mathcal{J}$\ generated in the crossed-product $C^{\ast }\left( 
\mathbb{F}_{2}\right) \rtimes _{\sigma }\mathbb{Z}_{2}$ by the spectral
projection $p=\frac{1}{2}\left( 1+W\right) $ of the canonical unitary $W$ in 
$C^{\ast }\left( \mathbb{F}_{2}\right) \rtimes \mathbb{Z}_{2}$ such that $%
WUW=V$. We first provide a simple yet useful description of $C^{\ast }\left( 
\mathbb{F}_{2}\right) \rtimes \mathbb{Z}_{2}$ in term of unitary generators,
before providing a description of the ideal $\mathcal{J}$ which will ease
the calculation of its $K$-theory.

\begin{lemma}
\label{ZZ2}Let $\sigma $ be the automorphism defined by $\sigma (U)=V$ and $%
\sigma (V)=U$ on the universal C*-algebra $C^{\ast }\left( \mathbb{F}
_{2}\right) $ generated by two universal unitaries $U$ and $V$. By
definition, $C^{\ast }\left( \mathbb{F}_{2}\right) \rtimes _{\sigma }\mathbb{%
\ Z}_{2}$ is the universal C*-algebra generated by three unitaries $U,V$ and 
$W $ subjects to the relations $W^{2}=1$ and $WUW^{\ast }=V$. The
C*-crossed-product $C^{\ast }\left( \mathbb{F}_{2}\right) \rtimes _{\sigma } 
\mathbb{Z}_{2}$ is equal to $C^{\ast }(U,W)$, or equivalently is the
universal C*-algebra generated by two unitaries $U$ and $W$ with the
relation $W^{2}=1$, or equivalently $C^{\ast }\left( \mathbb{F}_{2}\right)
\rtimes _{\sigma }\mathbb{Z}_{2}$ is *-isomorphic to $C^{\ast }\left( 
\mathbb{Z}\ast \mathbb{Z}_{2}\right) $.
\end{lemma}

\begin{proof}
The C*-subalgebra $C^{\ast }(U,W)$ of $C^{\ast }\left( \mathbb{F}_{2}\right)
\rtimes _{\sigma }\mathbb{Z}_{2}$ contains $WUW=V$ and thus equals $C^{\ast
}\left( \mathbb{F}_{2}\right) \rtimes _{\sigma }\mathbb{Z}_{2}$.

Let us now prove that the C*-subalgebra $C^{\ast }(U,W)$ is universal for
the given relations. Let $u,w$ be two arbitrary unitaries in some arbitrary
C*-algebra such that $w^{2}=1$. Let $v=wuw\in C^{\ast }(u,w)$. By
universality of the crossed-product $C^{\ast }\left( \mathbb{F}_{2}\right)
\rtimes _{\sigma }\mathbb{Z}_{2}$ there exists a unique *-morphism $\varphi
:C^{\ast }\left( \mathbb{F}_{2}\right) \rtimes _{\sigma }\mathbb{Z}
_{2}\longrightarrow C^{\ast }(u,w)$ such that $\varphi (U)=u$, $\varphi
(V)=v $ and $\varphi (W)=w$. Thus $C^{\ast }(U,W)$ is universal for the
proposed relations. In particular, it is *-isomorphic (by uniqueness of the
universal C*-algebra for the given relations) to $C^{\ast }\left( \mathbb{Z}%
\ast \mathbb{Z}_{2}\right) $.
\end{proof}

Now, the ideal $\mathcal{J}$ can be described as the kernel of a
particularly explicit *-morphism of $C^{\ast }\left( \mathbb{F}_{2}\right)
\rtimes _{\sigma }\mathbb{Z}_{2}$.

\begin{proposition}
\label{MoritaEq}Let $\sigma $ be the automorphism defined by $\sigma (U)=V$
and $\sigma (V)=U$ on the universal C*-algebra $C^{\ast }\left( \mathbb{F}
_{2}\right) $ generated by two universal unitaries $U$ and $V$. Let $W$ be
the canonical unitary of the crossed-product $C^{\ast }\left( \mathbb{F}
_{2}\right) \rtimes _{\sigma }\mathbb{Z}_{2}$ such that $WUW=V$. The fixed
point C*-algebra $C^{\ast }\left( \mathbb{F}_{2}\right) _{1}$ is strongly
Morita equivalent to the kernel $\mathcal{J}$ in $C^{\ast }\left( \mathbb{F}
_{2}\right) \rtimes _{\sigma }\mathbb{Z}_{2}$ of the *-morphism $\varphi
:C^{\ast }\left( \mathbb{F}_{2}\right) \rtimes _{\sigma }\mathbb{Z}
_{2}\longrightarrow C(\mathbb{T})$ defined by $\varphi (W)=-1$ and $\varphi
(U)(z)=\varphi (V)(z)=z$ for all $z\in \mathbb{T}$.
\end{proposition}

\begin{proof}
Let $\varphi $ be the unique *-morphism from $C^{\ast }(U,W)$ into $C(%
\mathbb{T})$ defined using the universal property of Lemma (\ref{ZZ2}) by: $%
\varphi (W)(z)=-1$ and $\varphi (U)(z)=z$ for all $z\in \mathbb{T}$. Since $%
\varphi (p)=0$ by construction, $\mathcal{J}\subseteq \ker \varphi $. We
wish to show that $\ker \varphi \subseteq \mathcal{J}$ as well.

Let $a\in \ker \varphi $. Let $\pi $ be a representation of $C^{\ast }\left( 
\mathbb{F}_{2}\right) \rtimes _{\sigma }\mathbb{Z}_{2}$ which vanishes on $%
\mathcal{J}$. Then $\pi (p)=0$ so $\pi (W)=-1$. Hence, for any $b\in C^{\ast
}\left( \mathbb{F}_{2}\right) _{-1}$ (see Lemma(\ref{split})) then $\pi
(b)=\pi (W)\pi (b)\pi (W)=\pi (-b)$ so that $\pi (b)=0$. In particular, $\pi
(U-V)=0$ and thus $\pi (U)=\pi (V)$. Thus, up to unitary equivalence, $\pi
=s\circ \varphi $ where $s:C(\mathbb{T})\rightarrow C(S)$ with $S$ the
spectrum of $\pi (U)$ and $s$ is the canonical surjection. In particular, $%
\pi (a)=s\circ \varphi (a)=0$. Since $\pi $ is arbitrary, we conclude that
the image of $a$ in $C^{\ast }\left( \mathbb{F}_{2}\right) \rtimes _{\sigma }%
\mathbb{Z}_{2}/\mathcal{J}$ is null, and thus $a\in \mathcal{J}$ as required.
\end{proof}

Our goal now is to compute the $K$-theory of the fixed point C*-algebra $%
C^{\ast }\left( \mathbb{F}_{2}\right) _{1}$. Since $C^{\ast }\left( \mathbb{%
F }_{2}\right) _{1}$ and the ideal $\mathcal{J}$ are Morita equivalent, they
have the same $K$-theory. By Proposition (\ref{MoritaEq}), we have the short
exact sequence $0\rightarrow \mathcal{J}\rightarrow C^{\ast }\left( \mathbb{%
F }_{2}\right) \rtimes _{\sigma }\mathbb{Z}_{2}\overset{\varphi }{%
\rightarrow } C(\mathbb{T})\rightarrow 0$, and it seems quite reasonable to
use the six-term exact sequence of $K$-theory to deduce the $K$-groups of $%
\mathcal{J }$ from the $K$-groups of $C^{\ast }\left( \mathbb{F}_{2}\right)
\rtimes _{\sigma }\mathbb{Z}_{2}$, as long as the later can be computed. The
next section precisely follows this path, starting by computing the $K$%
-theory of the crossed-product $C^{\ast }\left( \mathbb{F}_{2}\right)
\rtimes _{\sigma } \mathbb{Z}_{2}$.

\section{$K$-theory of the C*-crossed product and the Fixed Point C*-algebra}

We use the homotopy-based result in \cite{Cuntz82} to compute the $K$-theory
of the crossed-product $C^{\ast }\left( \mathbb{F}_{2}\right) \rtimes
_{\sigma }\mathbb{Z}_{2}$.

\begin{proposition}
\label{Kcrossed}Let $C^{\ast }\left( \mathbb{F}_{2}\right) =C^{\ast }(U,V)$
be the universal C*-algebra generated by two unitaries $U$ and $V$ and let $%
\sigma $ be the order-2 automorphism of $C^{\ast }\left( \mathbb{F}
_{2}\right) $ defined by $\sigma (U)=V$ and $\sigma (V)=U$. Then $%
K_{0}\left( C^{\ast }\left( \mathbb{F}_{2}\right) \rtimes _{\sigma }\mathbb{%
Z }_{2}\right) =\mathbb{Z}^{2}$ is generated by the spectral projections of $%
W$ and $K_{1}\left( C^{\ast }\left( \mathbb{F}_{2}\right) \rtimes _{\sigma } 
\mathbb{Z}_{2}\right) =\mathbb{Z}$ is generated by $U$.
\end{proposition}

\begin{proof}
By Lemma (\ref{ZZ2}), the C*-crossed-product $C^{\ast }\left( \mathbb{F}%
_{2}\right) \rtimes _{\sigma }\mathbb{Z}_{2}$ is *-isomorphic to 
\begin{equation*}
C^{\ast }\left( \mathbb{Z}\ast \mathbb{Z}_{2}\right) =C^{\ast }\left( 
\mathbb{Z}\right) \ast _{\mathbb{C}}C^{\ast }\left( \mathbb{Z}^{2}\right) =C(%
\mathbb{T})\ast _{\mathbb{C}}\mathbb{C}^{2}
\end{equation*}
where the free product is amalgated over the C*-algebra generated by the
respective units in each C*-algebra. More precisely, we embed $\mathbb{C}$
via, respectively, $i_{1}:\lambda \in \mathbb{C}\mapsto \lambda 1\in C(%
\mathbb{T})$ and $i_{2}:\lambda \in \mathbb{C}\mapsto \left( \lambda
,\lambda \right) \in \mathbb{C}^{2}$. There are natural *-morphisms from $C(%
\mathbb{T})$ and from $\mathbb{C}^{2}$ onto $\mathbb{C}$ defined
respectively by $r_{1}:f\in C(\mathbb{T})\mapsto f(1)$ and $r_{2}:\lambda
\oplus \mu \mapsto \lambda $. Now, $r_{1}\circ i_{1}=r_{2}\circ i_{2}$ is
the identity on $\mathbb{C}$. By \cite{Cuntz82}, we conclude that the
following sequences for $\varepsilon =0,1$ are exact: 
\begin{equation*}
0\rightarrow \mathbb{C}\overset{j_{\varepsilon }}{\rightarrow }%
K_{\varepsilon }(C(\mathbb{T})\oplus \mathbb{C}^{2})\overset{k_{\varepsilon }%
}{\rightarrow }K_{\varepsilon }\left( C(\mathbb{T})\ast _{\mathbb{C}}\mathbb{%
\ C}^{2}\right) \rightarrow 0
\end{equation*}%
where $j_{\varepsilon }=K_{\varepsilon }(i_{1})\oplus K_{\varepsilon
}(-i_{2})$ and $k_{\varepsilon }=K_{\varepsilon }(k_{1}+k_{2})$ where $k_{1}$
is the canonical embedding of $C(\mathbb{T})$ into $C(\mathbb{T})\ast _{%
\mathbb{C}}\mathbb{C}^{2}$ and $k_{2}$ is the canonical embedding of $%
\mathbb{C}^{2}$ into $C(\mathbb{T})\ast _{\mathbb{C}}\mathbb{C}^{2}$.

Now, $K_{1}(C(\mathbb{T})\oplus \mathbb{C}^{2})=\mathbb{Z}$ generated by the
identity $z\in \mathbb{T}\mapsto z\in C(\mathbb{T})$. Since $K_{1}\left( 
\mathbb{C}\right) =0$ we conclude that $K_{1}\left( C(\mathbb{T})\ast _{ 
\mathbb{C}}\mathbb{C}^{2}\right) =\mathbb{Z}$ generated by the canonical
unitary generator of $C(\mathbb{T})$. On the other hand, $K_{0}\left( C( 
\mathbb{T})\oplus \mathbb{C}^{2}\right) =\mathbb{Z}^{3}$ (where the first
copy of $\mathbb{Z}$ is generated by the unit $1_{C(\mathbb{T})}$ of $C( 
\mathbb{T})$ and the two other copies are generated by each of the
projections $(1,0)$ and $(0,1)$ in $\mathbb{C}^{2}$). Now, $\limfunc{ran}
j_{\varepsilon }$ is the subgroup of $\mathbb{Z}^{3}$ generated by the class
of $1_{C(\mathbb{T})}\oplus -1_{\mathbb{C}^{2}}$ where $1_{\mathbb{C}^{2}}$
is the unit of $\mathbb{C}^{2}$. This class is $(1,-1,-1)$, so we conclude
easily that $K_{0}\left( C(\mathbb{T})\oplus \mathbb{C}^{2}\right) =\mathbb{%
Z }^{2}$ is generated by the two projections in $\mathbb{C}^{2}$ (whose
classes are $(0,1,0)$ and $(0,0,1)$).

Using the *-isomorphism between $C(\mathbb{T})\oplus \mathbb{C}^{2}$ and $%
C^{\ast }\left( \mathbb{F}_{2}\right) \rtimes _{\sigma }\mathbb{Z}_{2}$ we
conclude that $K_{0}\left( C^{\ast }\left( \mathbb{F}_{2}\right) \rtimes
_{\sigma }\mathbb{Z}_{2}\right) =\mathbb{Z}^{2}$ is generated by the
spectral projections of $W$ while $K_{1}\left( C^{\ast }\left( \mathbb{F}
_{2}\right) \rtimes _{\sigma }\mathbb{Z}_{2}\right) =\mathbb{Z}$ is
generated by the class of $U$ (or $V$ as these are equal by construction).
\end{proof}

\begin{remark}
It is interesting to compare our results to the K-theory of the
C*-crossed-product by $\mathbb{Z}$ instead of $\mathbb{Z}_{2}$. One could
proceed with the standard six-terms exact sequence \cite[Theorem 10.2.1]%
{Blackadar98}, but it is even simpler to observe the following similar
result to Lemma (\ref{ZZ2}): the C*-crossed-product $C^{\ast }\left( \mathbb{%
\ F}_{2}\right) \rtimes _{\sigma }\mathbb{Z}$ is *-isomorphic to $C^{\ast } 
\mathbb{(F}_{2})$. If $C^{\ast }(\mathbb{F}_{2})$ is generated by the two
universal unitaries $U,V$ and $W$ is the canonical unitary in $C^{\ast
}\left( \mathbb{F}_{2}\right) \rtimes _{\sigma }\mathbb{Z}$ such that $WUW=V$
then once again $C^{\ast }(U,V,W)=C^{\ast }(U,W)$ and, following a similar
argument as for Lemma\ (\ref{ZZ2}) we observe that $C^{\ast }(U,W)$ is the
C*-algebra universal for two arbitrary unitaries.

Hence, $K_{0}\left( C^{\ast }\left( \mathbb{F}_{2}\right) \rtimes _{\sigma } 
\mathbb{Z}\right) =\mathbb{Z}$ is generated by the identity while $%
K_{1}\left( C^{\ast }\left( \mathbb{F}_{2}\right) \rtimes _{\sigma }\mathbb{%
Z }\right) =\mathbb{Z}^{2}$ is generated by $U$ and $W$.
\end{remark}

We can now deduce the $K$-theory of $C^{\ast }\left( \mathbb{F}_{2}\right)
_{1}$ from Proposition (\ref{MoritaEq}) and Theorem\ (\ref{Kcrossed}). The
natural way to do so is by using the six-terms exact sequence in $K$-theory
applied to the short exact sequence $\mathcal{J}\hookrightarrow C^{\ast
}\left( \mathbb{F}_{2}\right) \rtimes _{\sigma }\mathbb{Z}_{2}\overset{
\varphi }{\twoheadrightarrow }C(\mathbb{T})$.

\begin{theorem}
\label{KMorita}Let $C^{\ast }\left( \mathbb{F}_{2}\right) =C^{\ast }(U,V)$
be the universal C*-algebra generated by two unitaries $U$ and $V$ and let $%
\sigma $ be the order-2 automorphism of $C^{\ast }\left( \mathbb{F}
_{2}\right) $ defined by $\sigma (U)=V$ and $\sigma (V)=U$. Let $C^{\ast
}\left( \mathbb{F}_{2}\right) _{1}=\left\{ a\in C^{\ast }\left( \mathbb{F}
_{2}\right) :\sigma (a)=a\right\} $ be the fixed point C*-algebra for $%
\sigma $. Then $K_{0}(C^{\ast }\left( \mathbb{F}_{2}\right) _{1})=\mathbb{Z}$
is generated by the identity in $C^{\ast }\left( \mathbb{F}_{2}\right) _{1}$
and $K_{1}\left( C^{\ast }\left( \mathbb{F}_{2}\right) _{1}\right) =0$.
\end{theorem}

\begin{proof}
By \cite[Proposition 3.4]{Rieffel80}, $C^{\ast }\left( \mathbb{F}_{2}\right)
_{1}$ is strongly Morita equivalent to the ideal $\mathcal{J}$\ generated in 
$C^{\ast }\left( \mathbb{F}_{2}\right) \rtimes _{\sigma }\mathbb{Z}_{2}$ by
the projection $p=\frac{1}{2}\left( 1+W\right) $. Using Proposition (\ref%
{MoritaEq}), we can apply the six-terms exact sequence to the short exact
sequence 
\begin{equation*}
0\longrightarrow \mathcal{J}\overset{i}{\longrightarrow }C^{\ast }\left( 
\mathbb{F}_{2}\right) \rtimes _{\sigma }\mathbb{Z}_{2}\overset{\varphi }{
\longrightarrow }C(\mathbb{T})\longrightarrow 0
\end{equation*}
where $i$ is the canonical injection and $\varphi $ is the *-morphism of
Proposition (\ref{MoritaEq}).

We denote by $\Phi $ the quotient isomorphism from $C^{\ast }\left( \mathbb{F%
}_{2}\right) \rtimes _{\sigma }\mathbb{Z}_{2}/\mathcal{J}$ onto $C(\mathbb{T}%
)$ induced by $\varphi $. Since we know the $K$-theory of $C^{\ast }\left( 
\mathbb{F}_{2}\right) \rtimes _{\sigma }\mathbb{Z}_{2}$ by Theorem (\ref%
{Kcrossed}), including a set of generators of the $K$-groups, and the $K$
-theory of $C(\mathbb{T})$, we can easily deduce the $K$-theory of $\mathcal{%
\ J}$. Indeed, we have the following exact sequence \cite[9.3 p. 67]%
{Blackadar98}: 
\begin{equation}
\begin{array}{ccccc}
K_{0}\left( \mathcal{J}\right) & \overset{K_{0}(i)}{\longrightarrow } & 
K_{0}\left( C^{\ast }\left( \mathbb{F}_{2}\right) \rtimes _{\sigma }\mathbb{Z%
}_{2}\right) =\mathbb{Z}^{2} & \overset{K_{0}(\varphi )}{\longrightarrow } & 
K_{0}\left( C\left( \mathbb{T}\right) \right) \\ 
\delta \uparrow &  &  &  & \downarrow \beta \\ 
K_{1}\left( C\left( \mathbb{T}\right) \right) & \overset{K_{1}(\varphi )}{%
\longleftarrow } & K_{1}\left( C^{\ast }\left( \mathbb{F}_{2}\right) \rtimes
_{\sigma }\mathbb{Z}_{2}\right) =\mathbb{Z} & \overset{K_{1}(i)}{%
\longleftarrow } & K_{1}\left( \mathcal{J}\right)%
\end{array}
\label{exact}
\end{equation}%
Each statement in the following argument follows from the exactness of (\ref%
{exact}).Trivially, $K_{0}\left( \varphi \right) $ is a surjection, so $%
\beta =0$. Hence $K_{1}(i)$ is injective. Yet, as $\varphi (U):z\in \mathbb{%
\ T}\mapsto z$, we conclude that $K_{1}\left( \varphi \right) $ is an
isomorphism (since it maps a generator to a generator), and thus $K_{1}(i)=0$
. Now $K_{1}(i)=0$ and $\beta =0$ implies that $K_{1}\left( \mathcal{J}%
\right) =0$. Since $K_{1}(\varphi )$ is surjective, $\delta =0$ and thus $%
K_{0}(i)$ is injective. Its image is thus isomorphic to $K_{0}\left( 
\mathcal{J}\right) $ and coincide with $\ker K_{0}\left( \varphi \right) $.
Now, $K_{0}(\varphi )(p)=0$ and $K_{0}(\varphi )(1-p)=1$ (by Theorem\ (\ref%
{Kcrossed}), $p$ and $1-p$ generate $K_{0}\left( C^{\ast }\left( \mathbb{F}%
_{2}\right) \rtimes _{\sigma }\mathbb{Z}_{2}\right) $). Hence the image of $%
K_{0}(i)$ is isomorphic to the copy of $\mathbb{Z}$ generated by $1-p$ in $%
K_{0}\left( C^{\ast }\left( \mathbb{F}_{2}\right) \rtimes _{\sigma }\mathbb{Z%
}_{2}\right) $.
\end{proof}

We can go a little deeper in the structure of the fixed point C*-algebra $%
C^{\ast }\left( \mathbb{F}_{2}\right) _{1}$. Of interest is to compare the
ideal $\mathcal{J}$ and its natural restrictions to $C^{\ast }\left( \mathbb{%
\ F}_{2}\right) $ and $C^{\ast }\left( \mathbb{F}_{2}\right) _{1}$. The
motivation for this comparison is to understand the obstruction to the
existence of any nontrivial unitary in $C^{\ast }\left( \mathbb{F}
_{2}\right) _{1}$ in the sense of $K$-theory.

\begin{theorem}
\label{K0II1}Let $\theta $ be the *-epimorphism $C^{\ast
}(U,V)\twoheadrightarrow C(\mathbb{T})$ defined by $\varphi (U)=\varphi
(V):z\in \mathbb{T}\mapsto z$. Let $\mathcal{I}=\ker \theta $. Then $K_{0}(%
\mathcal{I})=0$ and $K_{1}\left( \mathcal{I}\right) =\mathbb{Z}$ where the
generating unitary in $\mathcal{I}^{+}=\mathcal{I}+1$ of the $K_{1}$ group
is $UV^{\ast }$.

Let $\mathcal{I}_{1}=C^{\ast }\left( \mathbb{F}_{2}\right) _{1}\cap \mathcal{%
\ I}$. We then have $C^{\ast }\left( \mathbb{F}_{2}\right) _{1}/\mathcal{I}
_{1}=C(\mathbb{T})$. Then $K_{1}\left( \mathcal{I}_{1}\right) =0$ while $%
K_{0}\left( \mathcal{I}_{1}\right) =\mathbb{Z}$.
\end{theorem}

\begin{proof}
We first calculate the $K$-theory of $\mathcal{I}$. This can be achieved in
at least two natural ways: by means of exact sequences (see Remark (\ref%
{IexactSeq})) or directly, by the following simple argument. Let $%
K_{1}\left( i\right) :K_{1}\left( \mathcal{I}\right) \longrightarrow
K_{1}\left( C^{\ast }\left( \mathbb{F}_{2}\right) \right) $ be the $K_{1}$
-lift of the canonical inclusion $i:\mathcal{I}\longrightarrow C^{\ast
}\left( \mathbb{F}_{2}\right) $. We first will identify the range of $%
K_{1}(i)$ and then show that $K_{1}(i)$ is injective.

Let $Z\in M_{n}\left( \mathcal{I}\right) +1_{n}$ be a unitary for some $n\in 
\mathbb{N}$. In particular, $Z$ is a unitary in $M_{n}\left( C^{\ast }\left( 
\mathbb{F}_{2}\right) \right) $. Let $(k,k^{\prime })\in \mathbb{Z}
^{2}=K_{1}\left( C^{\ast }\left( \mathbb{F}_{2}\right) \right) $ be the
class $[Z]_{C^{\ast }\left( \mathbb{F}_{2}\right) }$ of $Z$ in $K_{1}\left(
C^{\ast }\left( \mathbb{F}_{2}\right) \right) $. Now, by definition of $%
\theta $ and $\mathcal{I}$, we have $K_{1}\left( \theta \right)
([Z]_{C^{\ast }\left( \mathbb{F}_{2}\right) })=K_{1}\left( [1]_{C\left( 
\mathbb{T}\right) }\right) =0$. Yet, since $\theta (U)=\theta (V)=z\in 
\mathbb{T}\mapsto z$, we conclude that $K_{1}\left( \theta \right) \left(
k,k^{\prime }\right) =k+k^{\prime }$. Hence $[Z]_{C^{\ast }\left( \mathbb{F}
_{2}\right) }=\left( k,-k\right) $ for some $k\in \mathbb{Z}$. Moreover, $%
UV^{\ast }-1\in \mathcal{I}$ by construction, and $[UV^{\ast }]_{C^{\ast
}\left( \mathbb{F}_{2}\right) }=(1,-1)$, so the range of $K_{1}(i)$ is the
subgroup generated by $(1,-1)$ in $\mathbb{Z}^{2}=K_{1}\left( C^{\ast
}\left( \mathbb{F}_{2}\right) \right) $.

On the other hand, assume now that $[Z]_{C^{\ast }\left( \mathbb{F}
_{2}\right) }=0$, i.e. $Z\oplus 1_{m-n}$ is connected to $1$ in $M_{m}\left(
C^{\ast }\left( \mathbb{F}_{2}\right) \right) $ for some $m\geq n$. To ease
notation, let $1_{k}$ be the identity in $M_{k}\left( C^{\ast }\left( 
\mathbb{F}_{2}\right) \right) $, let $Y_{1}=Z\oplus 1_{m-n}$ and let $\left(
Y_{t}\right) _{t\in \lbrack 0,1]}$ be the homotopy of unitaries joining $%
Y_{1}$ to $Y_{0}=1_{m}$ in $M_{m}\left( C^{\ast }\left( \mathbb{F}
_{2}\right) \right) $. Let $\Xi $ be the the unique *-endomorphism of $%
C^{\ast }\left( \mathbb{F}_{2}\right) $ defined by $\Xi (U)=\Xi (V)=U$. For
any *-endomorphism $\eta $ of $C^{\ast }\left( \mathbb{F}_{2}\right) $ we
let $M_{m}\left( \eta \right) $ be the canonical *-endomorphism of $%
M_{m}\left( C^{\ast }\left( \mathbb{F}_{2}\right) \right) $ induced from $%
\eta $. Set $Y_{t}^{\prime }=Y_{t}(M_{n}(\Xi )(Y_{t}))^{\ast }$ for all $%
t\in \lbrack 0,1]$. Then we check immediately that $M_{n}(\theta
)(Y_{t})=1_{m}$, so $Y_{t}-1_{m}\in M_{m}\left( \mathcal{I}\right) $ by
definition. Hence, $\left( Y_{t}^{\prime }\right) _{t\in \lbrack 0,1]}$ is
now an homotopy in $M_{m}\left( \mathcal{I}\right) +1_{m}$ between $%
Y_{1}^{\prime }$ and $Y_{0}^{\prime }$. Since $M_{n}(\Xi )(Y_{1})=1_{m}$ by
construction since $Y_{1}\in C^{\ast }\left( \mathbb{F}_{2}\right) _{1}$, we
conclude that $Y_{1}^{\prime }=Y_{1}$ and $Y_{0}=Y_{0}^{\prime }$ (as $\Xi
(1)=1$). Hence, $[Z]_{\mathcal{I}}=0$. Therefore, $K_{1}(i)$ is injective.
Consequently, $K_{1}\left( \mathcal{I}\right) =\mathbb{Z}$ generated by the
class of $UV^{\ast }$.

A similar argument applies to homotopy of projections. We deduce that $%
K_{0}\left( \mathcal{I}\right) =\left\{ 0\right\} $.

Now, we turn to the ideal $\mathcal{I}_{1}$. We recall from Lemma (\ref%
{split}) the notation $C^{\ast }\left( \mathbb{F} _{2}\right) _{-1}=\left\{
a-\sigma (a):a\in C^{\ast }\left( \mathbb{F} _{2}\right) \right\} $. Our
first observation is that $C^{\ast }\left( \mathbb{F}_{2}\right)
_{-1}\subseteq \mathcal{I}$. Indeed, since $\theta (U)=\theta (V)$ we have $%
\theta \circ \sigma =\theta $ and thus $\theta (a-\sigma (a))=0$ for all $%
a\in C^{\ast }\left( \mathbb{F}_{2}\right) $. Therefore, $C^{\ast }\left( 
\mathbb{F}_{2}\right) /\mathcal{I}=C^{\ast }\left( \mathbb{F}_{2}\right)
_{1}/\left( \mathcal{I}\cap C^{\ast }\left( \mathbb{F}_{2}\right)
_{1}\right) $ since $C^{\ast }\left( \mathbb{F}_{2}\right) =C^{\ast }\left( 
\mathbb{F}_{2}\right) _{1}\oplus C^{\ast }\left( \mathbb{F}_{2}\right) _{-1}$
as vector spaces by Lemma (\ref{split}).

Using the six-terms exact sequence and Theorem (\ref{KMorita}), we can
compute the $K$-theory of the ideal $\mathcal{I}_{1}$: 
\begin{equation*}
\begin{array}{ccccc}
K_{0}\left( \mathcal{I}_{1}\right) & \overset{K_{0}(i)}{\longrightarrow } & 
K_{0}\left( C^{\ast }\left( \mathbb{F}_{2}\right) _{1}\right) =\mathbb{Z} & 
\overset{K_{0}(q)}{\longrightarrow } & K_{0}\left( C\left( \mathbb{T}\right)
\right) =\mathbb{Z} \\ 
\delta \uparrow &  &  &  & \downarrow \beta \\ 
K_{1}\left( C\left( \mathbb{T}\right) \right) =\mathbb{Z} & \overset{%
K_{1}\left( q\right) }{\longleftarrow } & K_{1}\left( C^{\ast }\left( 
\mathbb{F}_{2}\right) _{1}\right) =0 & \overset{K_{1}\left( i\right) }{%
\longleftarrow } & K_{1}\left( \mathcal{I}_{1}\right)%
\end{array}%
\end{equation*}%
where $i$ and $q$ are again the canonical injection and surjection. Each
subsequent argument follows from the exactness of the six-terms sequence.
Since $K_{0}\left( C^{\ast }\left( \mathbb{F}_{2}\right) _{1}\right) $ is
generated by the class of the unit in $C^{\ast }\left( \mathbb{F}_{2}\right)
_{1}$ and $q(1)=1$ generated $K_{0}\left( C(\mathbb{T})\right) $, we
conclude that $K_{0}(q)$ is the identity, so $K_{0}(i)=0$ and $\beta =0$.
Hence $K_{0}\left( \mathcal{I}_{1}\right) =\delta (\mathbb{Z})$ and $%
K_{1}(i) $ is injective. Since $K_{1}\left( C^{\ast }\left( \mathbb{F}%
_{2}\right) _{1}\right) =0$ and $\beta =0$ we conclude that $K_{1}\left( 
\mathcal{I}_{1}\right) =0$. Therefore $K_{1}(q)=0$, so $\delta $ is
injective and we get $K_{0}\left( \mathcal{I}_{1}\right) =\mathbb{Z}$.
\end{proof}

\begin{remark}
\label{IexactSeq}There is an alternative calculation of the $K$-theory of
the ideal $\mathcal{I}$ using the simple six-term exact sequence: 
\begin{equation*}
\begin{array}{ccccc}
K_{0}(\mathcal{I})=0 & \overset{K_{0}(i)}{\longrightarrow } & K_{0}\left(
C^{\ast }\left( \mathbb{F}_{2}\right) \right) =\mathbb{Z} & \overset{K_{0}(q)%
}{\longrightarrow } & K_{0}\left( C(\mathbb{T})\right) =\mathbb{Z} \\ 
\delta =0\uparrow &  &  &  & \downarrow \beta =0 \\ 
K_{1}\left( C(\mathbb{T})\right) =\mathbb{Z} & \overset{K_{1}(q)}{%
\longleftarrow } & K_{1}\left( C^{\ast }\left( \mathbb{F}_{2}\right) \right)
=\mathbb{Z}^{2} & \overset{K_{1}(i)}{\longleftarrow } & K_{1}\left( \mathcal{%
I}\right) =\mathbb{Z}%
\end{array}%
\end{equation*}%
corresponding to the defining exact sequence $\mathcal{I}\overset{i}{%
\hookrightarrow }C^{\ast }\left( \mathbb{F}_{2}\right) \overset{q}{%
\twoheadrightarrow }C(\mathbb{T})$ with $i$ the canonical injection and $q$
the canonical surjection. Now, $K_{0}\left( C^{\ast }\left( \mathbb{F}%
_{2}\right) \right) $ is generated by $1$, and as $q(1)=1$ we see that $%
K_{0}(q)$ is the identity. Hence $\beta =0$ and $K_{0}(i)=0$ so $K_{0}\left( 
\mathcal{I}\right) =\delta (\mathbb{Z})$. On the other hand, $K_{1}\left(
C^{\ast }\left( \mathbb{F}_{2}\right) \right) $ is generated by $U$ and $V$,
respectively identified with $(1,0)$ and $(0,1)$ in $\mathbb{Z}^{2}$. We
have $q(U)=q(V):z\mapsto z$ which is the generator of $K_{1}\left( C(\mathbb{%
T})\right) $, so $K_{1}(q)$ is surjective and thus $\delta =0$. Hence $%
K_{0}\left( \mathcal{I}\right) =0$. On the other hand, $\ker K_{1}(q)$ is
the group generated by $(1,-1)$, the class of $UV^{\ast }$. Thus $K_{1}(i)$,
which is an injection, is in fact a bijection from $K_{1}\left( \mathcal{I}%
\right) $ onto its range $\ker K_{1}(q)$ and thus $K_{1}\left( \mathcal{I}%
\right) =\mathbb{Z}$ generated by $UV^{\ast }$, as indeed $q(UV^{\ast }-1)=0$
and thus $UV^{\ast }-1\in \mathcal{I}$.
\end{remark}

It is not too surprising that $K_{1}\left( \mathcal{I}_{1}\right) =0$ since $%
K_{1}\left( \mathcal{I}\right) =\mathbb{Z}$ is generated by $UV^{\ast }(-1)$
which is an element in $C^{\ast }\left( \mathbb{F}_{2}\right) _{-1}$ of
class $(1,-1)$ in $K_{1}\left( C^{\ast }\left( \mathbb{F}_{2}\right) \right) 
$, so it is not connected to any unitary in $C^{\ast }\left( \mathbb{F}
_{2}\right) _{1}$. Of course, this is not a direct proof of this fact, as
homotopies in $M_{n}(\mathcal{I}_{1})+I_{n}$ ($n\in \mathbb{N}$) is a more
restrictive notion than in $M_{n}\left( C^{\ast }\left( \mathbb{F}
_{2}\right) _{1}\right) $ ($n\in \mathbb{N}$). But more remarkable is the
fact that $K_{0}(\mathcal{I}_{1})$ contains some nontrivial element. Of
course, $\mathcal{I}_{1}$ is projectionless since $C^{\ast }\left( \mathbb{F}
_{2}\right) $ is by \cite{ChoiFree}, so the projection generating $K_{0}( 
\mathcal{I}_{1})$ is at least (and in fact, exactly in) $M_{2}\left( 
\mathcal{I}_{1}\right) +I_{2}$. We now turn to an explicit description of
the generator of $K_{0}(\mathcal{I}_{1})$ and we investigate why this
projection is trivial in both $K_{0}\left( C^{\ast }\left( \mathbb{F}
_{2}\right) _{1}\right) $ and $K_{0}\left( \mathcal{I}\right) $ but not in $%
K_{0}\left( \mathcal{I}_{1}\right) $. By exactness of the six-terms exact
sequence, this projection is exactly the obstruction to the nontriviality of 
$K_{1}\left( C^{\ast }\left( \mathbb{F}_{2}\right) _{1}\right) $.

\begin{theorem}
\label{Kgenerator}Let $\theta :C^{\ast }\left( \mathbb{F}_{2}\right)
\longrightarrow C(\mathbb{T})$ be the *-homomorphism defined by $\theta
(U)=\theta (V):z\in \mathbb{T}\mapsto z$. Let $\mathcal{I}_{1}=\left\{ a\in
C^{\ast }\left( \mathbb{F}_{2}\right) _{1}:\theta (a)=0\right\} $. Let $Z=%
\frac{1}{2}\left( U+V\right) \in C^{\ast }\left( \mathbb{F}_{2}\right) _{1}$%
. Let $\beta $ be the projection in $M_{2}\left( \mathcal{I}_{1}\right)
+I_{2}$ defined by: 
\begin{equation*}
\beta =\left[ 
\begin{array}{cc}
Z^{\ast }Z & Z^{\ast }\left( \sqrt[2]{1-ZZ^{\ast }}\right) \\ 
\left( \sqrt[2]{1-ZZ^{\ast }}\right) Z & 1-ZZ^{\ast }%
\end{array}%
\right] \text{.}
\end{equation*}%
Then the generator of $K_{0}\left( \mathcal{I}_{1}\right) $ is $[\beta ]_{%
\mathcal{I}_{1}}-[p_{2}]_{\mathcal{I}_{1}}$ where $p_{2}=\left[ 
\begin{array}{cc}
1 & 0 \\ 
0 & 0%
\end{array}%
\right] $ and, for any C*-algebra $A$, we denote by $[q]_{A}$ the $K_{0}$
-class of any projection $q\in M_{n}(A)$ for any $n\in \mathbb{N}$. On the
other hand, the projection $\beta $ is homotopic to $p_{2}$ in $M_{2}(%
\mathcal{I})+I_{2}$ and in $M_{2}(C^{\ast }\left( \mathbb{F}_{2}\right)
_{1}) $. Thus $[\beta ]_{\mathcal{I}}=0$ in $K_{0}(\mathcal{I})$ and $[\beta
]_{C^{\ast }\left( \mathbb{F}_{2}\right) _{1}}=[p_{2}]_{C^{\ast }\left( 
\mathbb{F}_{2}\right) _{1}}=[1]_{C^{\ast }\left( \mathbb{F}_{2}\right) _{1}}$
in $K_{0}\left( C^{\ast }\left( \mathbb{F}_{2}\right) _{1}\right) $.
\end{theorem}

A simple calculation shows that $\beta $ is a projection. We organize the
proof of Theorem (\ref{Kgenerator}) in several lemmas. We start with the two
quick observations that $\beta $ is homotopic to $p_{2}$ in $M_{2}(C^{\ast
}\left( \mathbb{F}_{2}\right) _{1})$ and in $M_{2}(\mathcal{I})+I_{2}$, and
then we prove that $[\beta ]_{\mathcal{I}_{1}}\not=[p_{2}]_{\mathcal{I}_{1}}$
in $K_{0}\left( \mathcal{I}_{1}\right) $.

\begin{lemma}
\label{Kgenerator0}The projection $\beta $ and the projection $1-p_{2}$ are
homotopic in $M_{2}(C^{\ast }\left( \mathbb{F}_{2}\right) _{1})$. Thus $%
[\beta ]_{C^{\ast }\left( \mathbb{F}_{2}\right) _{1}}=[1]_{C^{\ast }\left( 
\mathbb{F}_{2}\right) _{1}}$.
\end{lemma}

\begin{proof}
For all $t\in \lbrack 0,1]$ we set:%
\begin{equation*}
\beta _{t}=\left[ 
\begin{array}{cc}
t^{2}Z^{\ast }Z & tZ^{\ast }\sqrt[2]{1-tZZ^{\ast }} \\ 
\left( \sqrt[2]{1-tZZ^{\ast }}\right) tZ & 1-tZZ^{\ast }%
\end{array}%
\right] \text{.}
\end{equation*}%
Then $\left( \beta _{t}\right) _{t\in \lbrack 0,1]}$ is by construction an
homotopy in $M_{2}(C^{\ast }\left( \mathbb{F}_{2}\right) _{1})$ between $%
\beta _{1}=\beta $ and $\beta _{0}=1-p_{2}$. Trivially $1-p_{2}$ and $p_{2}$
are homotopic, and $[p_{2}]_{C^{\ast }\left( \mathbb{F}_{2}\right)
_{1}}=[1]_{C^{\ast }\left( \mathbb{F}_{2}\right) _{1}}$, hence our result.
\end{proof}

The important observation in the proof of Lemma (\ref{Kgenerator0}) is that
although $\beta _{t}\in M_{2}(C^{\ast }\left( \mathbb{F}_{2}\right) _{1})$
for all $t\in \lbrack 0,1]$, we have 
\begin{equation*}
\theta (tZ\sqrt[2]{1-tZZ^{\ast }})(z)=tz\sqrt[2]{1-t}\not=0
\end{equation*}
for $t\in (0,1)$ and $z\in \mathbb{T}$, so $\beta _{t}$ does not belong to
the ideals $\mathcal{I}$ and $\mathcal{I}_{1}$.

Since $K_{0}\left( \mathcal{I}\right) =0$ it is trivial that the class of $%
\beta $ in $K_{0}\left( \mathcal{I}\right) $ is null, but the exact reason
why it is so is interesting as a way to contrast with the calculations of
the class of $\beta $ in $K_{0}\left( \mathcal{I}_{1}\right) $.

\begin{lemma}
\label{Kgenerator3}In $M_{2}\left( \mathcal{I}\right) +I_{2}$ the projection 
$\beta $ is homotopic to $p_{2}$. Hence in $K_{0}\left( \mathcal{I}\right) $
we verify that we have indeed $[\beta ]_{\mathcal{I}}=0$.
\end{lemma}

\begin{proof}
The unitary equivalence in Lemma (\ref{Kgenerator2}) does not carry to the
unitalization of the ideal $\mathcal{I}$, but we can check that $\beta $ is
homotopic to $p_{2}$ in $M_{2}\left( \mathcal{I}\right) +I_{2}$. Set $Z_{t}=%
\frac{1}{2}\left( tU+(1-t)V\right) $ and set: 
\begin{equation*}
\beta _{t}=\left[ 
\begin{array}{cc}
Z_{t}^{\ast }Z_{t} & Z_{t}^{\ast }\left( \sqrt[2]{1-Z_{t}Z_{t}^{\ast }}%
\right) \\ 
\left( \sqrt[2]{1-Z_{t}Z_{t}^{\ast }}\right) Z_{t} & 1-Z_{t}Z_{t}^{\ast }%
\end{array}%
\right]
\end{equation*}%
for all $t\in \lbrack 0,1]$. As before, $\beta _{t}$ is a projection for all 
$t\in \lbrack 0,1]$ since $\left\Vert Z_{t}\right\Vert =1$ for all $t\in
\lbrack 0,1]$. Now, $\beta _{0}=\beta _{1}=\left[ 
\begin{array}{cc}
1 & 0 \\ 
0 & 0%
\end{array}%
\right] $ while $\beta _{\frac{1}{2}}=\beta $. Of course, $t\in \lbrack 0,%
\frac{1}{2}]\mapsto \beta _{t}$ is continuous. Moreover: 
\begin{equation*}
\left( tU+(1-t)V\right) \left( tU+(1-t)V\right) ^{\ast }=1+(t-t^{2})\left(
UV^{\ast }+VU^{\ast }\right)
\end{equation*}%
so $\theta (1-Z_{t}Z_{t}^{\ast })=0$. Hence, $\beta _{t}\in M_{2}\left( 
\mathcal{I}\right) +I_{2}$ for all $t\in \lbrack 0,1]$. Hence $\beta $ is
homotopic to $p_{2}$ in $M_{2}\left( \mathcal{I}\right) +I_{2}$.
\end{proof}

Unlike in the case of Lemma (\ref{Kgenerator0}), the homotopy used in the
proof of Lemma (\ref{Kgenerator2}) is in $M_{2}\left( \mathcal{I}\right) $,
but it is not in $M_{2}\left( C^{\ast }\left( \mathbb{F}_{2}\right)
_{1}\right) $ and hence not in $M_{2}\left( \mathcal{I}_{1}\right) $.

The crux of this matter is that $\beta $ is the obstruction to the existence
of a nontrivial element in $K_{1}\left( C^{\ast }\left( \mathbb{F}%
_{2}\right) _{1}\right) $. In view of Lemmas (\ref{Kgenerator0}) and (\ref%
{Kgenerator3}), we wish to see a concrete reason why $\beta $ can not have
the same class as $p_{2}$ in $K_{0}(\mathcal{I}_{1})$. We start with a
useful calculation: since $\beta $ and $p_{2}$ are homotopic in $C^{\ast
}\left( \mathbb{F}_{2}\right) _{1}$, they are unitarily equivalent as well,
and we now explicit a unitary implementing this equivalence:

\begin{lemma}
\label{Kgenerator2}Let $Y=\left[ 
\begin{array}{cc}
Z^{\ast } & \sqrt[2]{1-Z^{\ast }Z} \\ 
\sqrt[2]{1-ZZ^{\ast }} & -Z%
\end{array}%
\right] $. Then $Y$ is a unitary in $M_{2}\left( C^{\ast }\left( \mathbb{F}%
_{2}\right) _{1}\right) $ such that $Yp_{2}Y^{\ast }=\beta $.
\end{lemma}

\begin{proof}
Observe that $Z^{\ast }\left( 1-ZZ^{\ast }\right) =Z^{\ast }-Z^{\ast
}ZZ^{\ast }=(1-Z^{\ast }Z)Z^{\ast }$. Thus, for any $n\in \mathbb{N}$ we get
by a trivial induction that $Z^{\ast }\left( 1-ZZ^{\ast }\right)
^{n}=(1-Z^{\ast }Z)^{n}Z^{\ast }$. Hence, for any polynomial $p$ by
linearity, we have $Z^{\ast }\left( p\left( 1-ZZ^{\ast }\right) \right)
=\left( p(1-Z^{\ast }Z)\right) Z^{\ast }$. By Stone-Weierstrass, we deduce
that $Z^{\ast }f(1-ZZ^{\ast })=f(1-Z^{\ast }Z)Z^{\ast }$ for any continuous
function $f$ on the spectrum of $1-ZZ^{\ast }$ and $1-Z^{\ast }Z$ which is
the compact $[0,1]$, and in particular for the square root. Therefore:%
\begin{equation}
Z^{\ast }\sqrt[2]{1-ZZ^{\ast }}=\left( \sqrt[2]{1-Z^{\ast }Z}\right) Z^{\ast
}\text{.}  \label{eq1}
\end{equation}

Now, we have: 
\begin{eqnarray*}
YY^{\ast } &=&\left[ 
\begin{array}{cc}
Z^{\ast } & \sqrt[2]{1-Z^{\ast }Z} \\ 
\sqrt[2]{1-ZZ^{\ast }} & -Z%
\end{array}%
\right] \left[ 
\begin{array}{cc}
Z & \sqrt[2]{1-ZZ^{\ast }} \\ 
\sqrt[2]{1-Z^{\ast }Z} & -Z^{\ast }%
\end{array}%
\right] \\
&=&\left[ 
\begin{array}{cc}
Z^{\ast }Z+1-Z^{\ast }Z & 0 \\ 
0 & 1%
\end{array}%
\right]
\end{eqnarray*}%
using (\ref{eq1}) since $Z^{\ast }\sqrt[2]{1-ZZ^{\ast }}-\left( \sqrt[2]{%
1-Z^{\ast }Z}\right) Z^{\ast }=0$. Similarly, we get $Y^{\ast }Y=1_{2}$.

Now, we compute $Yp_{2}Y^{\ast }$: 
\begin{eqnarray*}
&&\left[ 
\begin{array}{cc}
Z^{\ast } & \sqrt[2]{1-Z^{\ast }Z} \\ 
\sqrt[2]{1-ZZ^{\ast }} & -Z%
\end{array}%
\right] \left[ 
\begin{array}{cc}
1 & 0 \\ 
0 & 0%
\end{array}%
\right] \left[ 
\begin{array}{cc}
Z & \sqrt[2]{1-ZZ^{\ast }} \\ 
\sqrt[2]{1-Z^{\ast }Z} & -Z^{\ast }%
\end{array}%
\right] \\
&=&\left[ 
\begin{array}{cc}
Z^{\ast } & \sqrt[2]{1-Z^{\ast }Z} \\ 
\sqrt[2]{1-ZZ^{\ast }} & -Z%
\end{array}%
\right] \left[ 
\begin{array}{cc}
Z & \sqrt[2]{1-ZZ^{\ast }} \\ 
0 & 0%
\end{array}%
\right] \\
&=&\left[ 
\begin{array}{cc}
Z^{\ast }Z & Z^{\ast }\sqrt[2]{1-ZZ^{\ast }} \\ 
\left( \sqrt[2]{1-ZZ^{\ast }}\right) Z & 1-ZZ^{\ast }%
\end{array}%
\right] =\beta \text{.}
\end{eqnarray*}

Last, we observe that $\sigma (Z)=Z$ by construction and thus $\sigma (Y)=Y$
as well: in other words, $Y\in M_{2}\left( C^{\ast }\left( \mathbb{F}%
_{2}\right) _{1}\right) $ (and we recover that $\beta $ is unitarily
equivalent in $M_{2}\left( C^{\ast }\left( \mathbb{F}_{2}\right) _{1}\right) 
$ to $p_{2}$).
\end{proof}

Note that $Z-\lambda 1\not\in \mathcal{I}$ for all $\lambda \in \mathbb{C}$
and so $Y$ does not belong to $M_{2}(\mathcal{I})+1_{2}$. Indeed, the
following lemma shows that $\beta $ and $p_{2}$ do not have the same $K$%
-class in $\mathcal{I}_{1}$, precisely because the conjunction of the
conditions of symmetry and being in the kernel of $\theta $ make it
impossible to deform one into the other, even though each condition alone
does not create any obstruction.

\begin{lemma}
\label{Kgenerator4}We have $[\beta ]_{\mathcal{I}_{1}}-[p_{2}]_{\mathcal{I}
_{1}}\not=0$ in $K_{0}\left( \mathcal{I}_{1}\right) $.
\end{lemma}

\begin{proof}
To prove Theorem\ (\ref{Kgenerator}), it remains to show that $[\beta ]_{%
\mathcal{I}_{1}}-[p_{2}]_{\mathcal{I}_{1}}$ is a generator for $K_{0}\left( 
\mathcal{I}_{1}\right) $. Let $\delta :K_{1}\left( C(\mathbb{T})\right)
\longrightarrow K_{0}\left( \mathcal{I}_{1}\right) $ be the exponential map
in the six-term exact sequence in $K$-theory induced by the exact sequence $%
0\rightarrow \mathcal{I}_{1}\rightarrow C^{\ast }\left( \mathbb{F}%
_{2}\right) _{1}\overset{\theta }{\rightarrow }C(\mathbb{T})\rightarrow 0$.
Let us denote by $z$ the canonical unitary $z:\omega \in \mathbb{T}\mapsto
\omega $ in $C(\mathbb{T})$. Let us also denote by $\theta _{2}$ the map
induced by $\theta $ on $M_{2}\left( C^{\ast }\left( \mathbb{F}_{2}\right)
\right) $. By \cite[Proposition 9.2.3]{Rordam00}, if $u$ is any unitary in $%
M_{2}(C^{\ast }\left( \mathbb{F}_{2}\right) _{1})$ such that $\theta
_{2}\left( u\right) =\left[ 
\begin{array}{cc}
z & 0 \\ 
0 & z^{\ast }%
\end{array}%
\right] $, then $\delta ([z]_{C(\mathbb{T})})=\left[ up_{2}u^{\ast }\right]
_{\mathcal{I}_{1}}-[p_{2}]_{\mathcal{I}_{1}}$. In particular, $\theta
_{2}(Y)=\left[ 
\begin{array}{cc}
z & 0 \\ 
0 & z^{\ast }%
\end{array}%
\right] $ so $\delta ([z]_{C(\mathbb{T})})=[Yp_{2}Y^{\ast }]_{\mathcal{I}%
_{1}}-[p_{2}]_{\mathcal{I}_{1}}=[\beta ]_{\mathcal{I}_{1}}-[p_{2}]_{\mathcal{%
I}_{1}}$.

On the other hand, by Theorem (\ref{K0II1}), $\delta $ is an isomorphism of
group. Since $[z]_{C(\mathbb{T})}$ is a generator of $K_{1}\left( C\left( 
\mathbb{T}\right) \right) $ we conclude that $[\beta ]_{\mathcal{I}%
_{1}}-[p_{2}]_{\mathcal{I}_{1}}$ is a generator of $K_{0}\left( \mathcal{I}%
_{1}\right) $.
\end{proof}

We thus have proven Theorem (\ref{Kgenerator}) by identifying $[\beta ]_{%
\mathcal{I}_{1}}-[p_{2}]_{\mathcal{I}_{1}}$ as the generator of $K_{0}\left( 
\mathcal{I}_{1}\right) $ and verifying that without the conjoint conditions
of symmetry via $\sigma $ and $\theta $, the difference of the classes of $%
\beta $ and $p_{2}$ is null in both $K_{0}\left( C^{\ast }\left( \mathbb{F}%
_{2}\right) _{1}\right) $ and in $K_{0}\left( \mathcal{I}\right) $.

\bibliographystyle{amsplain}
\bibliography{thesis}

\end{document}